\documentclass[12pt]{article}
\usepackage[a4paper, margin=2.5cm]{geometry}
\usepackage{amsmath,amssymb,amsfonts, color, inputenc}
\usepackage{authblk}
\usepackage{booktabs, url, bbm}

\usepackage{latexsym}
\usepackage{psfrag}
\usepackage{graphicx}
\usepackage[margin=1cm]{caption}

\newtheorem{theorem}{Theorem}[section]

\newtheorem{remark}[theorem]{Remark}

\title{Nonparametric estimation of marginal distributions for unordered pairs}
\author{Laura Dumitrescu and Darcy Harcourt \thanks{This paper is based
on the the second author's honours project.} \\
School of Mathematics and Statistics \\ Victoria University of Wellington, New Zealand}

\begin{document}
\maketitle

\abstract{In this article, we consider the estimation of the marginal distributions for pairs of data are recorded, with unobserved order in each pair. New estimators are proposed and their asymptotic properties are established, by proving a Glivenko-Cantelli theorem and a functional central limit result. Results from a simulation study are included and we illustrate the applicability of the method on the homologous chromosomes data.}

\noindent {\em Keywords}: asymptotic inference; empirical process; nonparametric estimation; unordered pairs.

\noindent {\em AMS Classification}: Primary 62G30; Secondary 62G20.

\section{Introduction}
For independent and identically distributed random variables $X_1, \ldots, X_n$ the classical approach to the estimation of the unknown distribution $F$ on $\mathbb{R}$ is to use the empirical distribution function, $F_n(x).$ This function estimates the arbitrary distribution function $F$ at a rate of $n^{-1/2}$ in the supnorm loss leading to statistical tests which are optimal, in the sense of having power against all alternatives that can be separated from the null distribution by at least a constant multiple of $n^{-1/2}$ (see e.g. \cite{gine-nickl16}).

Recently, a test based on a Kolmogorov-Smirnov type statistic was shown in \cite{dumitrescu-khmaladze19} to only distinguish alternatives with a separation boundary of order $n^{-1/4}$ when testing the equality of marginal distributions for bivariate data, with unobservable order in each pair. More precisely, given a sequence of paired observations $(X_i, Y_i)_{i \ge 1}$, with marginals $F_1$ and $F_2,$ respectively, when the order in each pair is not observed, it is no longer possible to evaluate the corresponding empirical distributions, $F_{1n}$ and $F_{2n},$ and an approach based on a symmetrized empirical process was proposed. Namely, the process $\mathbb{R}_n^s(u,v)$ was defined as the symmetrized version of $\mathbb{R}_n(x,y)={n}^{1/2}\{\mathbb{F}_n(x,y) - \mathbb{Q}_n(x,y)\},$ with $\mathbb{Q}_n= Q_n \times Q_n$ and $Q_n=(F_{1n} + F_{2n})/2.$ Its limiting distribution was obtained as a restriction of a Brownian pillow to the class of symmetric sets. The idea has further been used in \cite{roberts21} to propose the implementation of the Khmaladze rotation in higher dimensions. 

In this article, we show that estimating the marginal distributions of the bivariate data may be achieved under the restriction of observability and propose new estimators. We show that it is possible to clearly distinguish (and hence, estimate) between the two marginals, as long as one distribution dominates the other. However, as expected, the estimation becomes more difficult in the tail regions, as well as in any neighbouring region to the intersection of the marginals. We investigate the asymptotic properties of our proposed  estimators and obtain a Glivenko-Cantelli theorem, as well a functional central limit theorem result.  
 
Our research is motivated by applications on genetic data, particularly, in karyotype analysis, where measurements of different characteristics are collected on homologous chromosomes. In this case, the question of interest is to determine if there exist significant differences between the chromosomes derived from the mother and the chromosomes derived from the father (see \cite{BaChaBa17}) and if so, we propose the corresponding estimators for the marginal distributions of the paired measurements.  

Another example arises in controlled trials, where ``blinding'' is often employed. This refers to the situation where some or more individuals (the analyst and/or the participant) involved in the study are unaware of the assigned treatment (see e.g. \cite{DayAltman00} for a detailed description and motivation). Blind assessment of the treatment outcome during the course of a clinical trial comparing, for example, the effects of two treatments on individuals, generates unordered pairs of observations. The blind analysis of interim data is often required as, in some cases, there may be an increasing risk of further experimentation (see also \cite{BaChaBa17}, \cite{MiFrKi09}).


\section{Marginal distributions estimation}
\label{estimation}
To motivate the form of the proposed estimators, we begin with a few simple remarks. Let $(X, Y)$ be a bivariate random vector with independent components whose distribution functions are denoted as $F_1$ and $F_2,$ respectively. Denote by $F^{(1)}$ and $F^{(2)}$ the distribution functions of the random variables $\min\{X, Y\}$ and $\max\{X,Y\}$ so that 
\begin{equation*}
F^{(1)}(x) = F_1(x) + F_2(x) - F_1(x)F_2(x), \ F^{(2)}(x) = F_1(x) F_2(x), \  x \in \mathbb{R}.
\end{equation*}
The following inequalities are equivalent, for any $x \in \mathbb{R},$ 
\begin{equation}
F^{(2)}(x) \le \left\{ \frac{F^{(1)}(x)+F^{(2)}(x)}{2}\right\}^2, \ F_1(x) F_2(x) \le \left\{ \frac{F_{1}(x)+F_{2}(x)}{2}\right\}^2 \nonumber
\end{equation}
and equality holds if and only if $F_1=F_2.$

Based on $F^{(1)}$ and $F^{(2)},$ for arbitrary $x \in \mathbb{R}$, we consider the functions
\begin{eqnarray*}
G_{1,2}(x) &=& \frac{F^{(1)}(x) + F^{(2)}(x) \pm [\{F^{(1)}(x) + F^{(2)}(x)\}^2 - 4F^{(2)}(x)]^{1/2} }{2}  \\
&=& \frac{F_{1}(x) + F_{2}(x) \pm |F_1(x) - F_2(x)| }{2}. 
\end{eqnarray*} 
Their form shows that, unless one distribution dominates the other, the marginal distributions $F_1$ and $F_2$ cannot be uniquely determined from $F^{(1)}$ and $F^{(2)}.$ 

In the general case, for $x \in \mathbb{R}$, denote the minimum and maximum functions as $G_1(x)=\min\{F_1(x), F_2(x)\} $ and $G_2(x)=\max\{F_1(x), F_2(x)\}.$ It follows that $G_1(x)=\alpha(F^{(1)}(x), F^{(2)}(x))$ and $G_2(x)=\beta(F^{(1)}(x), F^{(2)}(x)),$ where 
$$\alpha(s,t)=\frac{s + t - \{(s+t)^2 - 4t\}^{1/2} }{2}, \ \beta(s,t)=\frac{s + t + \{(s+t)^2 - 4t\}^{1/2} }{2}$$ 
are defined on ${\cal D}=\{(s, t) \in [0,1]^2, \mbox{ such that } (s+t)^2 - 4t \ge 0\}.$   

We are interested in making inference based on a sample $\{(X_i,Y_i)\}_{1 \le i \le n}$ of paired measurements, with continuous marginals, and for which the order in each pair is not observed. Therefore, the data collected consists of a sequence of pairs $\{(U_{i}, V_{i})\}_{1 \le i \le n},$ where $U_{i}=\min\{X_i,Y_i\},$ $V_{i}=\max\{X_i,Y_i\}$ and the random variables $\{U_{i}\}_{1 \le i \le n}$ and $\{V_{i}\}_{1 \le i \le n}$ form two sequences of i.i.d. random variables. Let $\displaystyle{F^{(1)}_n(x)=\frac{1}{n}\sum_{i=1}^n \mathbbm{1}_{\{ U_i \le x\}}}$ and $\displaystyle{F^{(2)}_n(x)=\frac{1}{n}\sum_{i=1}^n \mathbbm{1}_{\{V_i \le x\}}}$ denote their empirical distribution functions, respectively. 

When the two marginal distributions are different, estimation will be possible on those intervals where there is a clear separation between the two. We denote ${\cal S} = \{x \in \mathbb{R}, \mbox { such that } \inf_{x \in \mathbb{R}}|F_1(x) - F_2(x)| \ge m, \mbox{where } m \in (0,1)\}.$

For each $x \in {\cal S},$ on the event $\Omega_{n}=\{\inf_{x \in {\cal S}}\{F^{(1)}_n(x) + F^{(2)}_n(x)\}^2 - 4F^{(2)}_n(x) \ge 0\},$ we consider the estimators $G_{1n}(x)=\alpha(F^{(1)}_n(x), F^{(2)}_n(x))$ and $G_{2n}(x)=\beta(F^{(1)}_n(x), F^{(2)}_n(x))$  
\begin{eqnarray*}
G_{1n}(x) &=& \frac{F^{(1)}_n(x) + F^{(2)}_n(x) - [\{F^{(1)}_n(x) + F^{(2)}_n(x)\}^2 - 4F^{(2)}_n(x)]^{1/2} }{2}, \\
G_{2n}(x) &=& \frac{F^{(1)}_n(x) + F^{(2)}_n(x) + [\{F^{(1)}_n(x) + F^{(2)}_n(x)\}^2 - 4F^{(2)}_n(x)]^{1/2} }{2}.
\end{eqnarray*}

\section{Consistency and limiting distribution}

In the first part of this section, we show that $G_{1n}(x)$ and $G_{2n}(x)$ are strongly consistent estimators of $G_1(x)$ and $G_2(x)$ and that the convergence is uniform. This result relies on the Glivenko-Cantelli theorem applied to $F^{(1)}_n$ and $F^{(2)}_n$, which gives the corresponding almost sure uniform convergence to $F^{(1)}$ and $F^{(2)},$ 
\begin{eqnarray}
\sup_{x \in \mathbb{R}} |F_{n}^{(1)}(x) - F^{(1)}(x)| \longrightarrow 0 , \ \sup_{x \in \mathbb{R}} |F_{n}^{(2)}(x) - F^{(2)}(x)|  \longrightarrow 0. \label{GC1}
\end{eqnarray}

\begin{theorem}
\label{theorem1}
The following properties are satisfied.
\begin{eqnarray*}
&& {\rm (i)} \mbox{ For large enough values of } n , \mbox{ the estimators } G_{1n}(x) \mbox{ and } G_{2n}(x) \mbox{ exist almost surely, }  x \in {\cal S}. \\
&& {\rm (ii)} \ \sup_{x \in {\cal S}} |G_{1n}(x) - G_1(x)|  \longrightarrow 0, \mbox{almost surely}, \\
&& {\rm (iii)} \ \sup_{x \in {\cal S}} |G_{2n}(x) - G_2(x)| \longrightarrow 0, \mbox{almost surely}. 
\end{eqnarray*}
\end{theorem}

{\bf Proof.}
(i). We show that, $pr(\Omega_n)=1,$ for large values of $n$. Let $0<\varepsilon< m^2$ be arbitrarily fixed and let ${\cal E}_1$ and ${\cal E}_2$ be the events on which the limits in \eqref{GC1} hold, respectively. Then, on ${\cal E}_1$, there exists $n_{1\varepsilon}$ such that 
\begin{equation}
\sup_{x \in \mathbb{R}} |F^{(1)}_n(x) - F^{(1)}(x)| < \varepsilon, \mbox{ for } n \ge n_{1\varepsilon} \label{eqF(1)}
\end{equation}
and on ${\cal E}_2$, there exists $n_{2\varepsilon}$ such that 
\begin{equation}
\sup_{x \in \mathbb{R}} |F^{(2)}_n(x) - F^{(2)}(x)| < \varepsilon, \mbox{ for } n \ge n_{2\varepsilon}. \label{eqF(2)}
\end{equation}

It follows that, on ${\cal E}_1 \cap {\cal E}_2,$ there exists $N_{\varepsilon}=\max \{n_{1\varepsilon}, n_{2\varepsilon}\}$ such that  \eqref{eqF(1)} and \eqref{eqF(2)} are satisfied, for every $n \ge N_{\varepsilon}.$ Since the function $\Delta(s,t)=(s+t)^2-4t$, defined on $[0,1]^2$, is uniformly continuous, there exist $\delta_{1\varepsilon}>0$ and $\delta_{2\varepsilon}>0$ such that, if $|s_n-s_0|<\delta_{1\varepsilon}$ and $|t_n-t_0|<\delta_{2\varepsilon},$ we have
$$ |\Delta(s_n, t_n) - \Delta(s_0, t_0)|<\varepsilon.$$

Let $x \in \mathbb{R}$ be arbitrary; taking $s_n=F^{(1)}_n(x),$ $s_0=F^{(1)}(x) ,$ $t_n=F^{(2)}_n(x)$ and $t_0=F^{(2)}(x)$, on the event $ {\cal E}_1 \cap {\cal E}_2,$ we have 
\begin{equation}
\sup_{x \in {\cal S}} \left|\{F^{(1)}_n(x) + F^{(2)}_n(x)\}^2 - 4F^{(2)}_n(x) - \{F_1(x) - F_2(x)\}^2  \right| < \varepsilon, \ n\ge N_{\varepsilon}, \label{ConvDet}
\end{equation}
so that $\{F^{(1)}_n(x) + F^{(2)}_n(x)\}^2 - 4F^{(2)}_n(x) > -\epsilon + m^2 > 0$ which concludes the proof of (i).

Since $\alpha$ and $\beta$ are continuous on the compact ${\cal D}=\{[2{t}^{1/2} - t, 1] \times [0,1], \ t\in [0,1]\},$ the proofs of (ii) and (iii) follow similarly.  \hfill $\Box$

\begin{remark}
{\rm (a) In the regions where $F_1$ and $F_2$ are closer, the estimator denoted as $D_n(x)=\{F^{(1)}_n(x) + F^{(2)}_n(x)\}^2 - 4F^{(2)}_n(x)$ would often have negative values, which we suggest to be truncated to 0, thus leading to approximately equal (to the average of between marginals) estimators $G_{1n}$ and $G_{2n}$. 

(b) Theorem \ref{theorem1} shows that estimators $G_{1n}$ and $G_{2n}$ exist almost surely and converge uniformly, on the set ${\cal S}$ where there is a clear discrepancy between the marginal distributions. However, on the intervals where $F_1$ and $F_2$ are close, is estimation still possible? We write $D_n(x),$ $x \in \mathbb{R}$ as follows 
\begin{eqnarray*}
	D_n(x) &=& \{F^{(1)}_n(x) + F^{(2)}_n(x) - F^{(1)}(x) - F^{(2)}(x)\}^2 \nonumber \\
	       &+& 2\{F^{(1)}(x) + F^{(2)}(x)\}\{F^{(1)}_n(x) + F^{(2)}_n(x) - F^{(1)}(x) - F^{(2)}(x)\}  \label{expansionDn}\\
	       &+& (-4)\{F^{(2)}_n(x) - F^{(2)}(x)\} +  \{F_{1}(x) - F_{2}(x)\}^2. \nonumber\\
	       &=& T_n(x) + \{F_{1}(x) - F_{2}(x)\}^2,
\end{eqnarray*}
where $T_n(x)$ denotes the sum of the first three terms of the expansion.
From Donsker's theorem (as in \eqref{D12}, below), 
$$\|F^{(1)}_n - F^{(1)}\|_{\infty} = O_{P}(n^{-1/2}), \ \|F^{(2)}_n - F^{(2)}\|_{\infty} = O_{P}(n^{-1/2})$$
so that $\|T_n\|_{\infty}=\|F_1+F_2\|_{\infty}O_{P}(n^{-1/2})$. Then, $D_n(x)$ is non-negative whenever $\{F_{1}(x) - F_{2}(x)\}^2$ is the dominant term: $$\displaystyle{\tilde{\cal S} =\{ x \in \mathbb{R} \mbox{ such that } \inf_{x \in \mathbb{R}} |F_1(x)-F_2(x)| \ge Cn^{-({1/4}-\delta)}, C>0, 0< \delta \le 1/4\}},$$ implying that $G_{1n}$ and $G_{2n}$ are consistent. This property is in line with a result in \cite{dumitrescu-khmaladze19} where it is shown that, under restricted observability of the order in pairs, the test statistic cannot distinguish alternatives converging to the null at rate $n^{-1/2}$ but with a rate which needs to decrease slower.}

\end{remark}

Next, we study the asymptotic distribution of the proposed estimators. With $x \in {\cal S}$, we introduce notations
\begin{eqnarray*}
h^{1-}(x) &=& 1 - \frac{F^{(1)}(x)+F^{(2)}(x)}{[\{F^{(1)}(x)+F^{(2)}(x)\}^2-4F^{(2)}(x)]^{1/2}}, \\
h^{2-}(x) &=& 1 - \frac{F^{(1)}(x)+F^{(2)}(x)-2}{[\{F^{(1)}(x)+F^{(2)}(x)\}^2-4F^{(2)}(x)]^{1/2}}, \\
h^{1+}(x) &=& 1 + \frac{F^{(1)}(x)+F^{(2)}(x)}{[\{F^{(1)}(x)+F^{(2)}(x)\}^2-4F^{(2)}(x)]^{1/2}}, \\
h^{2+}(x) &=& 1 + \frac{F^{(1)}(x)+F^{(2)}(x)-2}{[\{F^{(1)}(x)+F^{(2)}(x)\}^2-4F^{(2)}(x)]^{1/2}}.
\end{eqnarray*}

The first and second order partial derivatives of $\alpha,$ defined on ${\cal D}^*=\{(s,t) \in [0,1]^2 \mbox{ such that } (s+t)^2-4t > 0\}$, are as follows
\begin{eqnarray*}
\frac{\partial \alpha(s,t)}{\partial s} &=& \frac{1}{2} \left[ 1 - \frac{s+t}{\{(s+t)^2-4t\}^{1/2}}\right], \
\frac{\partial\alpha(s,t)}{\partial t} = \frac{1}{2} \left[ 1 - \frac{s+t-2}{\{(s+t)^2-4t\}^{1/2}}\right], \\
\frac{\partial^2\alpha(s,t)}{\partial s^2} &=& \frac{2t}{\{(s+t)^2-4t\}^{3/2}}, \
\frac{\partial^2\alpha(s,t)}{\partial t^2} = -\frac{2(s-1)}{\{(s+t)^2-4t\}^{3/2}}, \\
\frac{\partial^2\alpha(s,t)}{\partial s \partial t} &=& \frac{t-s}{\{(s+t)^2-4t\}^{3/2}}.
\end{eqnarray*}
 The standard result on empirical processes gives 
\begin{equation}
{n}^{1/2}(F_n^{(1)}-F^{(1)}) \Rightarrow v_{F^{(1)}}, \ {n}^{1/2}(F_n^{(2)}-F^{(2)}) \Rightarrow v_{F^{(2)}}, \mbox{ as } n \to \infty, \label{D12}
\end{equation}
on the space of c{\`a}dl{\`a}g functions, equipped with the Skorokhod topology. Here, $v_{F^{(1)}}(x)=B(F^{(1)}(x))$ and $v_{F^{(2)}}=B(F^{(2)}(x))$ are zero-mean Gaussian processes, where $B$ is a standard Brownian bridge on the unit interval. 

\begin{theorem}
As $n \to \infty,$ we have
\begin{eqnarray*}
{\rm (i)} \ {n}^{1/2}(G_{1n} - G_1) \Rightarrow \frac{1}{2} (v_{F^{(1)}} h^{1-} + v_{F^{(2)}} h^{2-}),\\
{\rm (ii)} \ {n}^{1/2}(G_{2n} - G_2) \Rightarrow \frac{1}{2} (v_{F^{(1)}} h^{1+} + v_{F^{(2)}} h^{2+}).\\
\end{eqnarray*}
\end{theorem}

{\bf Proof.} We only show part (i). Let $(s_0,t_0) \in {\cal D}^*;$ by Taylor's expansion of $\alpha$ around $(s_0, t_0)$, we write 
$$\alpha(s,t) = \alpha(s_0, t_0) + \begin{pmatrix} 
                                    s-s_0 \\ 
                                    t- t_0 
																	 \end{pmatrix}^T 
																	\bigtriangledown \alpha(s_0,t_0)
																	+ \frac{1}{2}\begin{pmatrix} s-s_0 \\
																	                             t- t_0 
																								\end{pmatrix}^T 
																		\bigtriangledown^2 \alpha(s^*,t^*)
																		\begin{pmatrix} 
																		s-s_0 \\ 
																		t- t_0 
																		\end{pmatrix},$$
where $|s^*-s_0| \le |s-s_0|$ and $|t^*-t_0| \le |t-t_0|.$ 


Let $x \in {\cal S}$ and take $s=F_n^{(1)}(x),$ $t=F_n^{(2)}(x),$ $s_0=F^{(1)}(x)$ and $t_0=F^{(2)}(x).$ For $\varepsilon >0,$ due to  \eqref{GC1}, there exists $n_{\varepsilon}$ such that, for every $n \ge n_{\varepsilon},$ we have $|s^*-F^{(1)}(x)| < \varepsilon$ and $|t^*-F^{(2)}(x)| < \varepsilon,$ almost surely, from where we obtain
\begin{eqnarray*}
 (s^*+t^*)^2-4t^* &<& \{F^{(1)}(x)+F^{(2)}(x) + 2 \varepsilon\}^2 - 4F^{(2)}(x) + 4\varepsilon \\
&=& \{F_1(x)-F_2(x)\}^2 + 4 \varepsilon^2 +4 \varepsilon\{F^{(1)}(x)+F^{(2)}(x)\} + 4\varepsilon \\
&\le &  1 + 4 \varepsilon^2+12\varepsilon,
\end{eqnarray*}
\begin{eqnarray*}
(s^*+t^*)^2-4t^* &>& \{F^{(1)}(x)+F^{(2)}(x) - 2 \varepsilon\}^2 - 4F^{(2)}(x) - 4\varepsilon \\
 &=& \{F_1(x)-F_2(x)\}^2 + 4 \varepsilon^2 - 4 \varepsilon\{F^{(1)}(x)+F^{(2)}(x)\} - 4\varepsilon \\
&\ge & m^2 + 4 \varepsilon^2 - 12\varepsilon.
\end{eqnarray*}
For $\varepsilon$ small enough such that $m^2 + 4 \varepsilon^2 - 12\varepsilon>0,$ i.e. $\varepsilon < \{3-(9-m^2)^{1/2}\}/2,$ it follows that the second order partial derivatives of $\alpha$ are uniformly bounded
\begin{eqnarray*}
\left|\frac{\partial^2\alpha(s^*,t^*)}{\partial s^2}\right|  &<&  \frac{2(1+\varepsilon)}{(m^2 + 4 \varepsilon^2 - 12\varepsilon)^{3/2}}, \ 
\left|\frac{\partial^2\alpha(s^*,t^*)}{\partial t^2}\right| <  \frac{2(1+\varepsilon)}{(m^2 + 4 \varepsilon^2 - 12\varepsilon)^{3/2}}, \\
\left|\frac{\partial^2\alpha(s^*,t^*)}{\partial s \partial t}\right| &<&  \frac{1+2\varepsilon}{(m^2 + 4 \varepsilon^2 - 12\varepsilon)^{3/2}}.
\end{eqnarray*}
This, together with \eqref{D12} justifies the uniformly negligibility of the second term of the right hand side of 
\begin{eqnarray*}
&&G_{1n}(x) - G_1(x) = \frac{1}{2}\begin{pmatrix} 
                                    F_n^{(1)}(x)-F^{(1)}(x) \\ 
                                    F_n^{(2)}(x)-F^{(2)}(x)
												\end{pmatrix}^T 
												\begin{pmatrix}
												              h^{1-}(x) \\
                                      h^{2-}(x)
												\end{pmatrix} \nonumber \\
									&&+			          
												\frac{1}{2}\begin{pmatrix} 
                                    F_n^{(1)}(x)-F^{(1)}(x) \\ 
                                    F_n^{(2)}(x)-F^{(2)}(x)
												\end{pmatrix}^T
												\bigtriangledown^2 \alpha(s^*,t^*)
												\begin{pmatrix} 
                                    F_n^{(1)}(x)-F^{(1)}(x) \\ 
                                    F_n^{(2)}(x)-F^{(2)}(x)
												\end{pmatrix}. \label{diff-exp} 
	\end{eqnarray*}
	 The conclusion now follows as an application of \eqref{D12}. 	\hfill $\Box$


\section{Simulated data}
\label{MCsim}
In this section we investigate the performance of the proposed estimators under two scenarios. Firstly, we consider the case when there is a dominance relation between the two marginals and take $F_1(x)=x$ and $F_2(x)=x^2.$   

We generate a sample of $n$ pairs $(X_i, Y_i),$ $i=1, \ldots, n$ where $X_i \sim F_1$ and $Y_i \sim F_2$ and calculate the empirical distribution functions $F_{1n},$ $F_{2n},$ $F_n^{(1)},$ $F_n^{(2)}$ of $X_i,$ $Y_i,$ $U_i$ and $V_i,$ as well as functions $G_{1n}$ and $G_{2n}$. The value of $\{F_n^{(1)}(x) + F_n^{(2)}(x)\}^2 - 4 F_n^{2}(x)$ was truncated to 0, whenever it was negative.

\begin{figure}[h!]
\includegraphics[width=\linewidth]{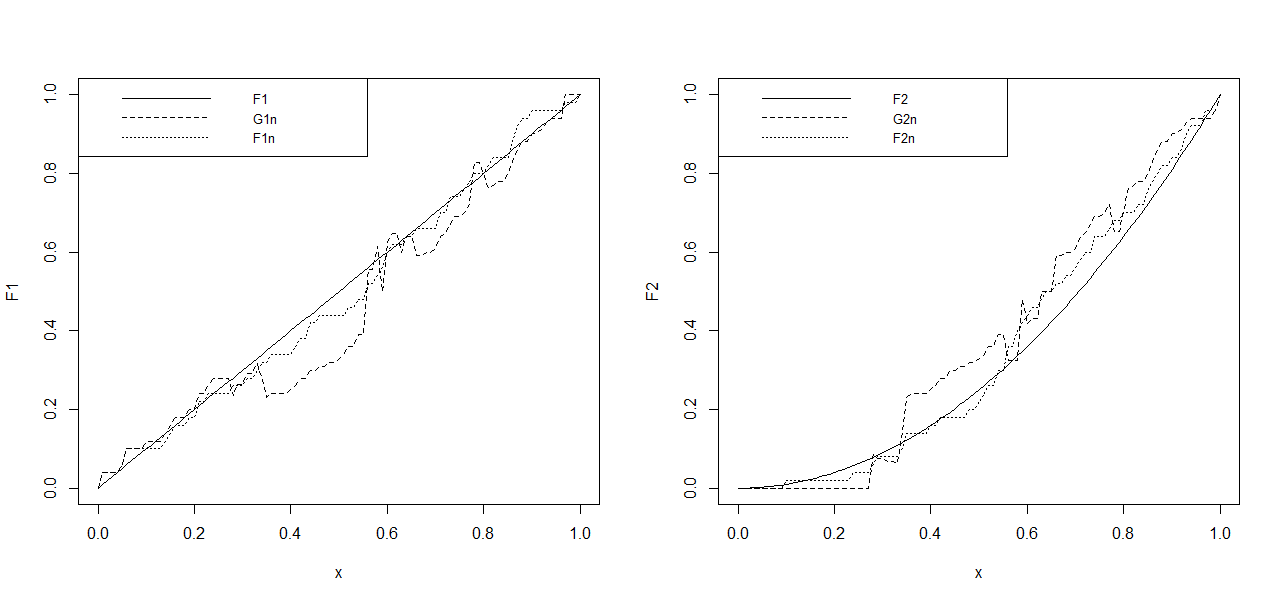}
\caption{Left panel: comparison between the graph of $F_1(x)=x$ (solid) and those of the estimators $G_{1n}$ (dashed) and $F_{1n}$ (dotted),  $n=50$. Right panel: comparison between the graph of $F_2(x)=x^2$ (solid) and those of the estimators $G_{2n}$ (dashed) and $F_{2n}$ (dotted), $n=50$.}
\label{fig:plotfns1}
\end{figure}

Figure \ref{fig:plotfns1} illustrates the graphs of the marginals and their estimators when $n=50$. The graphs of $F_{1n}$ and $F_{2n}$ are only added as a reference since they cannot be obtained when the order in each pair is not observed. As Figure \ref{fig:plotfns1} shows, the estimators $G_{1n}$ and $G_{2n}$ perform well in estimating the corresponding marginals $F_1$ and $F_2$, but not as well as $F_{1n}$ and $F_{2n}$, due to the restriction. Moreover, unlike the empirical distribution functions, the estimator functions $G_{1n}$ and $G_{2n}$ are not monotone.


\begin{figure}[h!]
\includegraphics[width=\linewidth]{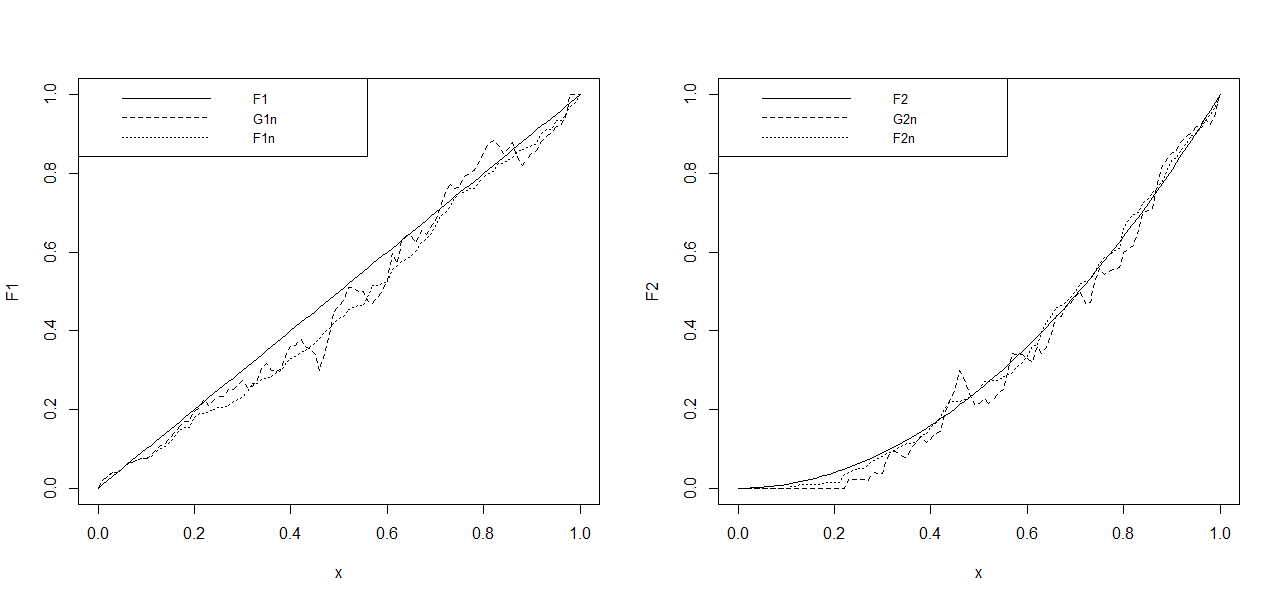}
\caption{Left panel: comparison between the graph of $F_1(x)=x$ (solid) and those of the estimators $G_{1n}$ (dashed) and $F_{1n}$ (dotted),  $n=200$. Right panel: comparison between the graph of $F_2(x)=x^2$ (solid) and those of the estimators $G_{2n}$ (dashed) and $F_{2n}$ (dotted), $n=200$.}
\label{fig:plotfns3}
\end{figure}

Figure \ref{fig:plotfns3} presents a similar conclusion, when $n=200$. As the sample size increases, the estimators $G_{1n}$ and $G_{2n}$ are shown to capture well the shape of $F_{1}$ and $F_2$ and further approaching them. 

As a reference, we added the empirical distributions of $U_i$ and $V_i,$ namely $F_n^{(1)}$ and $F_n^{(2)}$ to illustrate the fact that they are quite far away from the two marginal distributions and that, based only on this type of plot, a consistent estimation procedure is not evident. As remarked in \cite{dumitrescu-khmaladze19} (see e.g. their Figure 1), even when $F_1$ and $F_2$ are considerably different, obtaining estimators derived from only the distributions of the minima and maxima may be a challenging task.

\begin{figure}[h!]
\includegraphics[width=\linewidth]{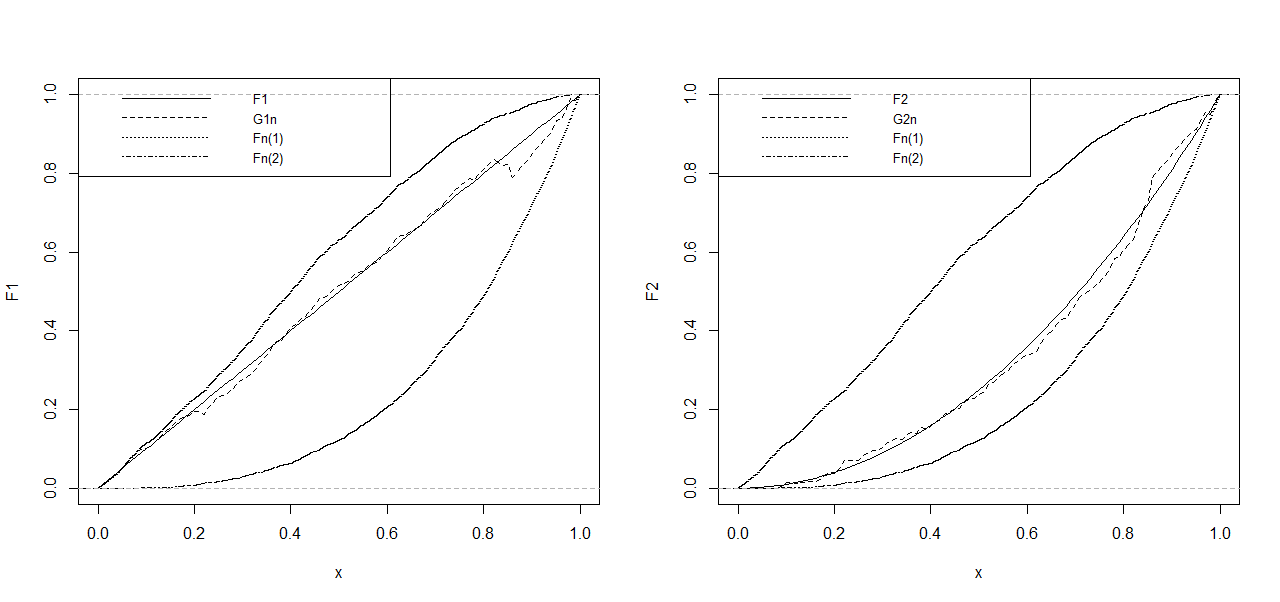}
\caption{Left panel: comparison between the graph of $F_1(x)=x$ (solid) and that of the estimator $G_{1n}$ (dashed), $n=2,000$. Right panel: comparison between $F_2(x)=x^2$ (solid) and that of the estimator $G_{2n}$ (dashed), $n=2,000$. In each panel, the graphs of $F^{(1)}_n$ (dotted) and $F_n^{(2)}$ (dotdashed) are added as a reference only.}
\label{fig:plotfns2}
\end{figure}

\begin{figure}[h!]
\includegraphics[width=\linewidth]{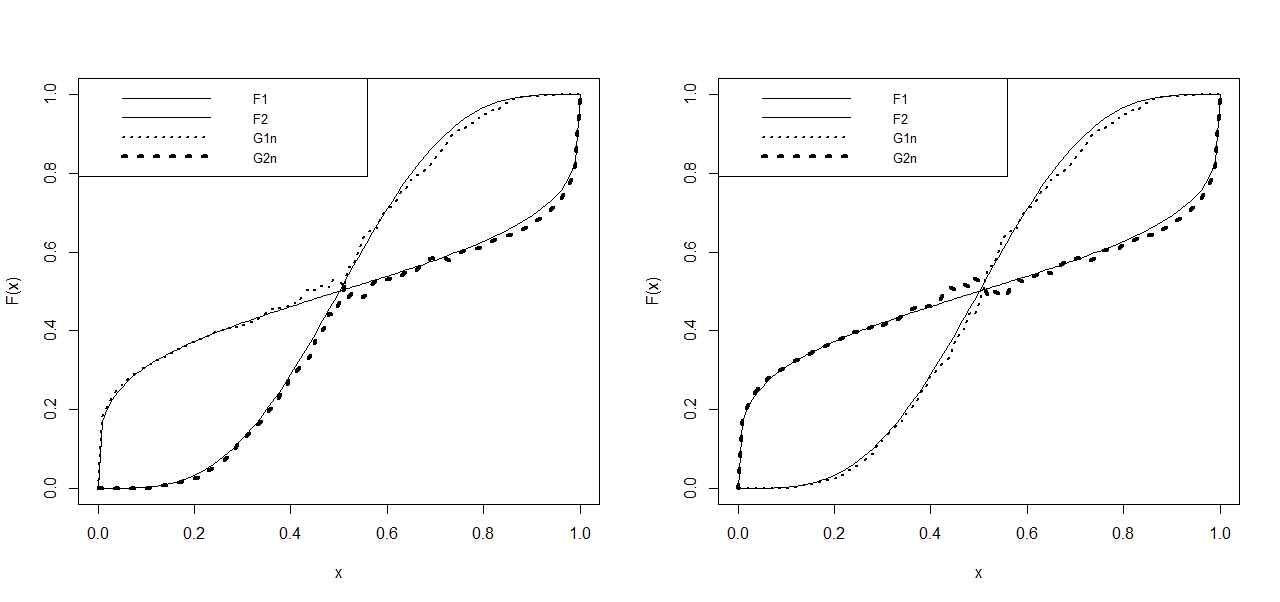}
\caption{Comparison between the graphs of the distribution function of $B(4,4)$ (denoted as $F_1$, solid) and the distribution function of $B(0.25,0.25)$ (denoted as $F_2$, solid) to those of the estimators $G_{1n}$ (dashed) and $G_{2n}$ (dashed), $n=2,000$.}
\label{fig:plotbetas}
\end{figure}

Secondly, we consider the case of two beta distributions and generate $n=2,000$ pairs of data $(X_i, Y_i)$, with $X_i \sim B(4,4)$ and $Y_i \sim B(0.25, 0.25)$. As Figure \ref{fig:plotbetas} shows, the estimators $G_{1n}$ and $G_{2n}$ follow closely the initial distributions even though, around the intersection point, distinguishing between the two marginals becomes more difficult. Moreover, in this case, estimation is not unique and two variants are possible.

\section{Application to the homologous chromosome dataset}
\label{realdata}
In each nucleus of a somatic human cell, there are 23 pairs of chromosomes and within each pair, one chromosome is derived from the mother DNA and the other is derived from the father DNA. Visually the chromosomes in the pair are not distinguishable and a sequence of unordered pairs of normalized measurements in the C-band area of the number 9 chromosome was analysed in \cite{li-qin11}. To determine if there exist significant differences between the chromosomes derived from the mother and those derived from the father, a test based on the empirical Shannon's mutual information was proposed. The data is included in Table 3 of \cite{li-qin11}. 

The approach proposed in \cite{li-qin11} is based on a semiparametric assumption on the probability densities of the marginal distributions, namely that the density ratio follows an exponential tilting model. The latter holds for normal densities with different means and same variance, exponential densities with different rates and the Poisson densities.

We consider a nonparametric test for $H_0:  F_1=F_2$ against $H_1:  F_1 \neq F_2$ and use the colour-blind empirical process proposed in \cite{dumitrescu-khmaladze19}
$$ \mathbb{R}_n^{s}(u,v)  = \mathbb{R}_n(u, v)+\mathbb{R}_n(v, u)-\mathbb{R}_n(u,u), \ u \le v ,$$
introduced as the restriction of   	
$$\mathbb{R}_n(x,y)={n}^{1/2} \left\{\mathbb{F}_n(x,y) - \frac{F_{1n} + F_{2n}}{2}(x)\frac{F_{1n} + F_{2n}}{2}(y) \right\}$$
to the class of symmetrized rectangles, where $\displaystyle{\mathbb{F}_{n}(x,y)= \frac{1}{n}\sum _{i=1}^{n}  \mathbbm{1}_{\{ X_{i} \le x, \ Y_{i} \le y\}}}.$

The limiting distribution of the process $\mathbb{R}_n^{s}$ corresponds to the that of a Brownian pillow on symmetric sets. A standard Brownian pillow is a Gaussian process $z(u,v)$ with covariance function $ \mathrm{E} \{z(u', v') z(u'', v'')\} = (\min \{u', u''\} - u'u'')(\min \{v', v''\} - v'v''), \ 0 \le u',u'',v',v'' \le 1$ and we approximate $\sup_{0 \le u\le v \le 1}{|\mathbb{R}_n^{s}(u,v)|}$ by $\sup_{0 \le u\le v \le 1}|z^{s}(u,v)|,$ where $z^{s}(u,v)=z(u,v)+z(v,u)-z(\min\{u,v\}, \min\{u,v\})$. Since the analytic form of the asymptotic distribution of the symmetrized Brownian pillow has not been derived, we propose the following simulation method (see also \cite{zhang14} for the standard Brownian pillow).

\begin{enumerate}
\item Generate $m^2$ independent and identically distributed random variables $\xi_{ij} \sim {\cal N}(0, 1/m^2),$ $1 \le i, j \le m.$ \\
\item Compute $\eta_{kl}=\sum_{i=1}^k \sum_{j=1}^l \xi_{ij},$ $1 \le k,l \le m.$ \\
\item Evaluate $\zeta_{kl}=\eta_{kl} - (l/m)\eta_{km} - (k/m) \eta_{ml} + (k/m)(l/m) \eta_{mm},$ $1 \le k,l \le m.$ \\
\item Take $\zeta^{s}_{kl}=\zeta_{kl}+\zeta_{lk} - \zeta_{kk} ,$ $1 \le k \le l \le m.$
\end{enumerate}

From the construction, $E(\zeta_{kl})=0$ and $cov(\zeta_{kl}, \zeta_{k'l'}) = (\min\{ k, k'\}/m - kk'/m^2)(\min\{l, l'\}/m - ll'/m^2)$ and so, for large values of $m,$ the distribution of $\sup_{1 \le k \le l \le m}|\zeta_{kl}^s|$ is approximately equal to that of $\sup_{0 \le u\le v \le 1}|z^{s}(u,v)|.$

The approximations of the $100(1-\alpha)\%$ upper quantiles (denoted as $z^{s}_{\alpha}$) of the distribution of $\sup_{0 \le u\le v \le 1}|z^{s}(u,v)|$ are given below.
\begin{center}
\begin{tabular}{c|c|c|ccc}
                 & $m$    & $n$     & $\alpha=0.1$ & $\alpha=0.05$ & $\alpha =0.01$ \\
\hline
$z^{s}_{\alpha}$ & 200    & 1,000    & 0.8592        & 0.9367         & 1.0489 \\
$z^{s}_{\alpha}$ & 300    & 1,000    & 0.8662        & 0.9390         & 1.0277 \\
$z^{s}_{\alpha}$ & 1,000  & 100,000  & 0.8868        & 0.9533         & 1.0804 \\
\end{tabular}.
\end{center}

For the homologous chromosomes data, the observed value of the test statistic $\sup_{u \le v} |\mathbb{R}_n^{s}(u,v)|$ is $1.04454$, which leads to a $p$-value of 0.01664 and hence, the null hypothesis of equality between the marginals is rejected. The result is in agreement with the one obtained in \cite{li-qin11}, where a $p$-value of 0.0141 was obtained based on the asymptotic distribution of their test statistic, whereas a $p$-value of 0.0225 was obtained using a bootstrap method.

\begin{figure}[h!]
\includegraphics[width=\linewidth]{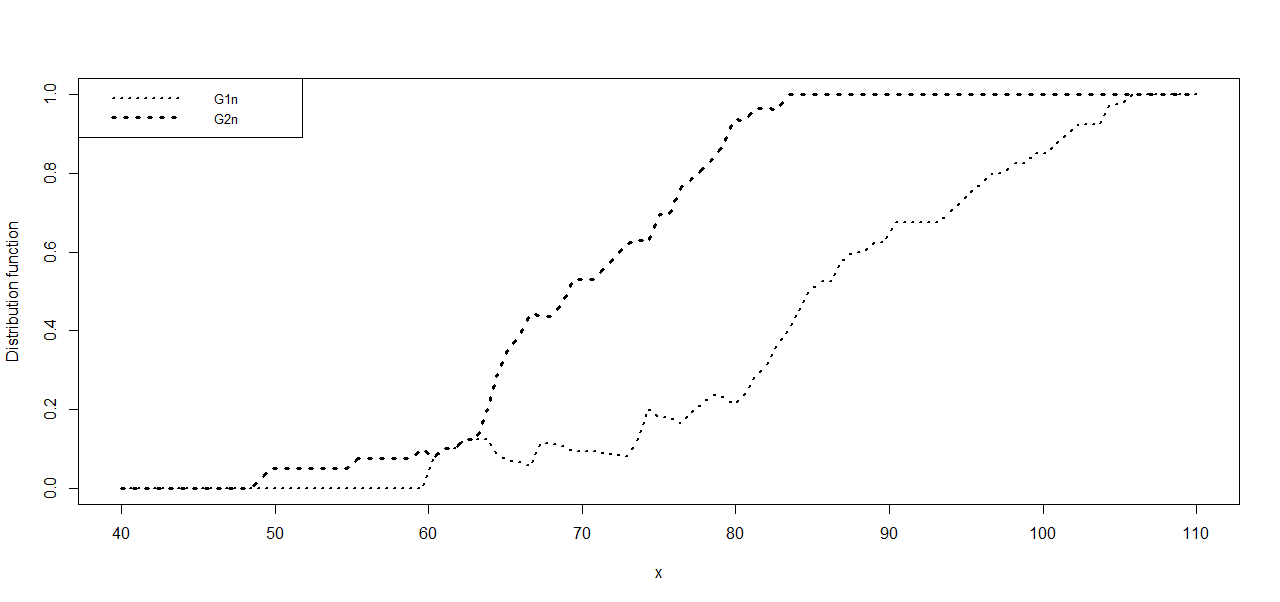}
\caption{Estimated marginal distribution functions for the homologous chromosome data. }
\label{fig:plotDF}
\end{figure}

With respect to the estimation, in Figure \ref{fig:plotDF} a plot of $G_{1n}$ and $G_{2n}$ is given, illustrating that, in the case of homologous chromosomes data, there is a clear separation between the two marginals. Moreover, it appears that one distribution dominates the other and that, as expected, discerning between the two becomes more difficult especially in the lower tail region.

\vspace{2mm}

\section*{Acknowledgement} 
The authors are thankful to Estate V. Khmaladze who suggested this problem and for related talks.


\begin{thebibliography}{7}
\expandafter\ifx\csname natexlab\endcsname\relax\def\natexlab#1{#1}\fi

\bibitem{BaChaBa17}
\textsc{Banerjee, T., Chattopadhyay, G. and Banerjee, K.} (2017). 
\newblock {Two stages test of means of unordered pairs.} 
\newblock {\em Statistics in Medicine}, {\bf 36}, 2466--2480.

\bibitem{DayAltman00} 
\textsc {Day, S. J. and Altman, D. G.} (2000). 
\newblock {Blinding in clinical trials and other studies.}
\newblock {\em BMJ}, {\bf 321}, 504.

\bibitem{dumitrescu-khmaladze19} 
\textsc {Dumitrescu, L. and Khmaladze, E. V.} (2019). 
\newblock {Asymptotic hypothesis testing for the colour blind problem.}
\newblock {\em Electronic Journal of Statistics,} {\bf 13}, 4573--4595.

\bibitem{gine-nickl16} 
\textsc {Giné, E. and Nickl, R.} (2016). 
\newblock {\em Mathematical foundations of infinite-dimensional statistical models}. 
\newblock Cambridge, University Press.

\bibitem{li-qin11} 
\textsc {Li, P. and Qin, J.} (2011). 
\newblock {A New Nuisance-Parameter Elimination Method With Application to the Unordered Homologous Chromosome Pairs Problem.} 
\newblock {\em Journal of the American Statistical Association,} {\bf 106}, 1476--1484. 

\bibitem{MiFrKi09} 
\textsc {Miller, F., Friede, T. and Kieser, M.} (2009). 
\newblock {Blinded assessment of treatment effects utilizing information about the randomization block length. }
\newblock {\em Statistics in Medicine}, {\bf 28}, 1690--1706.

\bibitem{roberts21} 
\textsc {Roberts, L. A.} (2021). 
\newblock {On the derivation of the Khmaladze transforms.} 
\newblock {\em arxiv}: 2101.07795.

\bibitem{zhang14} 
\textsc {Zhang, T.} (2014). 
\newblock {A Kolmogorov-Smirnov type test for independence between marks and points of marked points processes.}
\newblock {\em Electronic Journal of Statistics}, {\bf 8}, 2557--2584.
\end{thebibliography}
\end{document}